\documentclass{amsart}

\usepackage{amsmath}
\usepackage{amsthm}
\usepackage{amsfonts}
\usepackage{amssymb}
\usepackage{latexsym}
\usepackage{amscd}
\usepackage{stmaryrd}
\usepackage{mathrsfs}

\theoremstyle{plain}
\newtheorem{thm}{Theorem}[section]
\newtheorem{cor}[thm]{Corollary}
\newtheorem{lem}[thm]{Lemma}
\newtheorem{prop}[thm]{Proposition}

\def\proof{\noindent {\bf Proof.\;}}
\def\Li{\operatorname{Li}}

\newfont{\scyr}{wncyr10 scaled 550}
\def\shuffle{\,\mbox{\bf \scyr X}\,}

\numberwithin{equation}{section}

\begin{document}

\title{On harmonic sums and alternating Euler sums}

\date{}

\author{Zhong-hua Li}

\maketitle

\begin{center}
{\small Department of Mathematics, Tongji University, No. 1239 Siping Road},\\
{\small Shanghai 200092, China}\\
{\small Graduate School of Mathematical Sciences, The University of Tokyo}, \\
{\small 3-8-1 Komaba, Meguro, Tokyo 153-8914, Japan}\\
{\small E-mail address: lizhmath@gmail.com}
\end{center}

\vskip10pt

{\footnotesize
\begin{quote}

\noindent {\bf Abstract.} The explicit formulas expressing harmonic sums via alternating Euler sums (colored multiple zeta values) are given, and some explicit evaluations are given as applications.


\noindent{\bf 2010MSC}: 11M32, 40B05
\end{quote}
}


\section{Introduction}

Let $\mathbb{N}$ be the set of natural numbers. For $n,r\in\mathbb{N}$,
a generalized harmonic number $H_n^{(r)}$ is defined by
$$H_n^{(r)}:=\sum\limits_{j=1}^n\frac{1}{j^r},$$
which is a natural generalization of the harmonic number
$$H_n:=H_n^{(1)}=\sum\limits_{j=1}^n\frac{1}{j}.$$
Similarly, let
$$\overline{H}_n^{(r)}:=\sum\limits_{j=1}^n\frac{(-1)^{j+1}}{j^r}
$$
denote the alternating harmonic numbers.

Let $\mathbb{Z}$ be the set of integers. For $n, r_1,r_2,\ldots,r_l\in\mathbb{Z}\backslash\{0\}$ with
$n\geqslant 2$,
we define the harmonic sum
\begin{align}
S(n;r_1,r_2,\ldots,r_l):=\sum\limits_{k=1}^{\infty}\frac{X_k^{(r_1)}X_k^{(r_2)}\cdots X_k^{(r_l)}}{(k+1)^n},
\label{Eq:def-harmonicSum1}
\end{align}
where $X_k^{(r_i)}=H_k^{(r_i)}$ if $r_i>0$, and $X_k^{(r_i)}=\overline{H}_k^{(-r_i)}$ otherwise. In below, if $r<0$, we will denote it by $\overline{-r}$.
Similarly, for $n\in \mathbb{N}$, we define
\begin{align}
S(\bar{n};r_1,r_2,\ldots,r_l):=\sum\limits_{k=1}^{\infty}(-1)^{k+1}\frac{X_k^{(r_1)}X_k^{(r_2)}\cdots X_k^{(r_l)}}{(k+1)^n}.
\label{Eq:def-harmonicSum2}
\end{align}
We call $|n|+|r_1|+\cdots+|r_l|$ the weight of the harmonic sum $S(n;r_1,\ldots,r_l)$. Note that for any permutation $\sigma\in\mathfrak{S}_l$, we have
$$S(n;r_{\sigma(1)},\ldots,r_{\sigma(l)})=S(n;r_1,\ldots,r_l).$$

For example, for $n,m\in\mathbb{N}$ with $n\geqslant 2$, we have
$$S(n;\{1\}^m)=\sum\limits_{k=1}^\infty \left(1+\frac{1}{2}+\cdots+\frac{1}{k}\right)^m\frac{1}{(k+1)^n},$$
which is just $S_h(m,n)$ defined in \cite{BaileyBG,BBG}, and
$$S(n;\{\bar{1}\}^m)=\sum\limits_{k=1}^\infty \left(1-\frac{1}{2}+\cdots+\frac{(-1)^{k+1}}{k}\right)^m\frac{1}{(k+1)^n},$$
which is just $S_a(m,n)$ defined in \cite{BaileyBG,BBG}. For $n,m\in\mathbb{N}$, we have
$$S(\bar{n};\{1\}^m)
=\sum\limits_{k=1}^\infty \left(1+\frac{1}{2}+\cdots+\frac{1}{k}\right)^m\frac{(-1)^{k+1}}{(k+1)^n},$$
which is just $a_h(m,n)$ defined in \cite{BaileyBG,BBG}, and
$$S(\bar{n};\{\bar{1}\}^m)
=\sum\limits_{k=1}^\infty \left(1-\frac{1}{2}+\cdots+\frac{(-1)^{k+1}}{k}\right)^m\frac{(-1)^{k+1}}{(k+1)^n},$$
which is just $a_a(m,n)$ defined in \cite{BaileyBG,BBG}. Here we adopt the convenience that $\{s\}^m$ means the string $s$ repeats $m$ times.

The study of these harmonic sums was started by Euler. After that many different methods, including  partial fraction expansions, Eulerian Beta integrals, summation formulas for generalized hypergeometric functions and
contour integrals, have been used to evaluate these sums.
For details and historical introductions, please  see \cite{BaileyBG,BB,BBG,Choi-Srivastava,Chu,Doelder,Flajolet-Salvy,Shen} and references therein.

As in \cite{BorweinBradleyB}, for $k_1,\ldots,k_n\in \mathbb{N}$, we define the alternating Euler sums (colored multiple zeta values) by
\begin{align}
\zeta(k_1,\ldots,k_n;\sigma_1,\ldots,\sigma_n):=\sum\limits_{m_1>\cdots>m_n>0}
\frac{\sigma_1^{m_1}\cdots \sigma_n^{m_n}}{m_1^{k_1}\cdots m_n^{k_n}}, \;\;(k_1,\sigma_1)\neq(1,1),
\label{Eq:def-AEulerSum}
\end{align}
where $\sigma_i=\pm 1$ for all $1\leqslant i\leqslant n$. We call $k_1+\cdots+k_n$ the weight and $n$ the depth. As usual, if $\sigma_i=-1$
then $\bar{k}_i$ will be used. For example, we have
$$\zeta(\bar{1})=\zeta(1;-1)=\sum\limits_{m=1}^\infty \frac{(-1)^n}{n}=-\log 2.$$
If $\sigma_1=\cdots=\sigma_n=1$, then we get the well studied multiple zeta value
$$\zeta(k_1,\ldots,k_n)=\sum\limits_{m_1>\cdots>m_n>0}\frac{1}{m_1^{k_1}\cdots m_n^{k_n}},\;\;(k_1\geqslant 2).$$

When $l=1$, the harmonic sums $S(n;r_1)$ reduce to double alternating Euler sums.
In fact, we have $S(n;r)=\zeta(n,r)$ and
$S(n;\bar{r})=-\zeta(n,\bar{r})$ for $n\geqslant 2$, $r\geqslant 1$, and
$S(\bar{n};r)=\zeta(\bar{n},r)$ and $S(\bar{n};\bar{r})=-\zeta(\bar{n},\bar{r})$ for $n,r\geqslant 1$.
More generally, Flajolet and Salvy stated in \cite{Flajolet-Salvy} that
every $S(n;r_1,\ldots,r_l)$ for $n,r_1,\ldots,r_l\in\mathbb{N}$ with $n\geqslant 2$
is a $\mathbb{Q}$-linear combination of multiple zeta values.
However they didn't give the explicit formula and didn't use this fact as their main tools of \cite{Flajolet-Salvy}.
Now alternating Euler sums, especially multiple zeta values, are well developed since 1990s.
(One can find the list of references on multiple zeta values and Euler sums till now in the web of homepage
of Professor M. E. Hoffman.) So we think it is the time to give the explicit relation between the harmonic sums $S(n;r_1,\ldots,r_l)$ and
alternating Euler sums, and then treat $S(n;r_1,\ldots,r_l)$ using the results of alternating Euler sums.
This short note is devoted to this task. In Section 2, we give the explicit formulas expressing harmonic sums by alternating Euler sums.
In Section 3, we give some explicit evaluations as applications.


\section{Represented by alternating Euler sums}

Let $\mathscr{S}$ be the free abelian group generated by all finite sequences of nonzero integers. For
two sequences $\alpha=(k_1,\ldots,k_l),\beta=(k_{l+1},\ldots,k_{l+r})\in \mathscr{S}$, we define
$(\alpha,\beta)$ to be the sequence $(k_1,\ldots,k_l,k_{l+1},\ldots,k_{l+r})$ concatenating $\alpha$ and $\beta$.
Then for $\alpha=\sum a_i\alpha_i,
\beta=\sum b_j\beta_j\in\mathscr{S}$
with $a_i,b_j\in \mathbb{Z}$ and $\alpha_i,\beta_j$ being sequences, we define
$$(\alpha,\beta)=\sum a_ib_j (\alpha_i,\beta_j).$$
Now we define the harmonic shuffle product $\ast$ in $\mathscr{S}$
by linearities and the axioms:
\begin{itemize}
  \item $1\ast \alpha=\alpha\ast 1=\alpha$, for any sequence $\alpha$;
  \item $(k,\alpha)\ast(l,\beta)=(k,\alpha\ast (l,\beta))+(l,(k,\alpha)\ast\beta)+(p(k+l),\alpha\ast\beta)$, for
  any $k,l\in\mathbb{Z}\backslash\{0\}$ and any sequences $\alpha,\beta$,
\end{itemize}
where
$$p(k+l)=\left\{\begin{array}{ll}
k+l, & \text{\;if\;} k,l>0,\\
&\\
-k-l, & \text{\;if\;} k,l<0,\\
&\\
\overline{k-l}, & \text{\;if\;} k>0,l<0,\\
&\\
\overline{-k+l}, & \text{\;if\;} k<0, l>0,
\end{array}
\right.$$
and $1$ is the length zero sequence.
For example, we have
\begin{align*}
&(1)\ast(1)=(1,1)+(1,1)+(2)=2(1,1)+(2),\\
&(1)\ast (1)\ast (1)=2(1)\ast (1,1)+(1)\ast(2)=(3)+3(1,2)+3(2,1)+6(1,1,1),
\end{align*}
and
\begin{align*}
&(1)\ast(\bar{2})=(1,\bar{2})+(\bar{2},1)+(\bar{3}),\\
&(\bar{1})\ast (\bar{2})=(\bar{1},\bar{2})+(\bar{2},\bar{1})+(3).
\end{align*}
If we only consider sequences of positive integers, the harmonic product $\ast$ is just the one defined algebraically in \cite{Hoffman,IKZ}.

For any $\alpha=\sum a_i\alpha_i\in\mathscr{S}$, we define
$$\zeta(\alpha)=\sum a_i\zeta(\alpha_i)$$
if $\zeta(\alpha_i)$ exists for any $i$. For example,
if $\alpha=2(\bar{1},\bar{1})+(2)$, we have
$\zeta(\alpha)=2\zeta(\bar{1},\bar{1})+\zeta(2)$. Then it is easy to see that
$$\zeta(\alpha)\zeta(\beta)=\zeta(\alpha\ast\beta),$$
where $\alpha=(k_1,\ldots,k_n),\beta=(l_1,\ldots,l_m)\in\mathscr{S}$ with $k_1\geqslant 2$ if $k_1>0$ and
$l_1\geqslant 2$ if $l_1>0$.

Now by the definitions, we obtain the main theorem of this section.

\begin{thm}\label{Thm:rep-ast-general}
For any $n,r_1,r_2,\ldots,r_l\in\mathbb{Z}\backslash\{0\}$ with $n\geqslant 2$ if $n>0$, we have
\begin{align}
S(n;r_1,r_2,\ldots,r_l)=(-1)^k\zeta(n,(r_1)\ast(r_2)\ast\cdots \ast(r_l)),
\label{Eq:express-AEulerSum}
\end{align}
where $k=\sharp\{i\mid r_i<0\}$. In particular, every $S(n;r_1,\ldots,r_l)$ is a $\mathbb{Z}$-linear combination of alternating
Euler sums of weight $|n|+|r_1|+\cdots+|r_l|$ and depth at most $l+1$.
\end{thm}

A special case of the above theorem is

\begin{thm}\label{Thm:rep-ast}
For any $n,r_1,r_2,\ldots,r_l\in\mathbb{N}$ with $n\geqslant 2$, we have
\begin{align}
S(n;r_1,r_2,\ldots,r_l)=\zeta(n,(r_1)\ast(r_2)\ast\cdots \ast(r_l)).
\label{Eq:express-MZV}
\end{align}
In particular, every $S(n;r_1,\ldots,r_l)$ is a $\mathbb{Z}$-linear combination of multiple zeta values
of weight $n+r_1+\cdots+r_l$ and depth at most $l+1$.
\end{thm}

For example, we have
\begin{align}
S(n;\{1\}^2)=&\zeta(n,2)+2\zeta(n,1,1),\label{Eq:example-2}\\
S(n;\{1\}^3)=&\zeta(n,3)+3\zeta(n,2,1)+3\zeta(n,1,2)+6\zeta(n,1,1,1),\label{Eq:example-3}\\
S(n;\{\bar{1}\}^2)=&\zeta(n,2)+2\zeta(n,\bar{1},\bar{1}),\label{Eq:example-2-}\\
S(n;\{\bar{1}\}^3)=&-\zeta(n,\bar{3})-3\zeta(n,2,\bar{1})-3\zeta(n,\bar{1},2)-6\zeta(n,\bar{1},\bar{1},\bar{1}),\label{Eq:example-3-}
\end{align}
and
\begin{align*}
S(n;1,\bar{2})=-\zeta(n,\bar{3})-\zeta(n,1,\bar{2})-\zeta(n,\bar{2},1),\\
S(n;\bar{1},\bar{2})=\zeta(n,3)+\zeta(n,\bar{1},\bar{2})+\zeta(n,\bar{2},\bar{1}).
\end{align*}

In the special case $r_1=r_2=\ldots=r_l$, we can get more explicit formula as stated in the following corollary.

\begin{cor}\label{Cor:express-special}
For any $n\in\mathbb{Z}\backslash\{0\}$ and $k,r\in\mathbb{N}$ with $n\geqslant 2$ if $n>0$, we have
\begin{align}
S(n;\{r\}^k)
=&\sum\limits_{l=1}^k\sum\limits_{a_1+\cdots+a_l=k\atop a_i>0}\frac{k!}{a_1!\cdots a_l!}
\zeta(n,ra_1,ra_2,\ldots,ra_l),
\label{Eq:express-r+}
\end{align}
and
\begin{align}
S(n;\{\bar{r}\}^k)
=&(-1)^k\sum\limits_{l=1}^k\sum\limits_{a_1+\cdots+a_l=k\atop a_i>0}\frac{k!}{a_1!\cdots a_l!}
\zeta(n,\widetilde{ra_1},\widetilde{ra_2},\ldots,\widetilde{ra_l}),
\label{Eq:express-r-}
\end{align}
where we define
\begin{align}
\widetilde{ra_i}=\left\{\begin{array}{ll}
ra_i,&\text{\;if\;\,} a_i \text{\;is even},\\
&\\
\overline{ra_i},& \text{\;if\;\,} a_i \text{\;is odd}.
\end{array}\right.
\label{Eq:tilde}
\end{align}
\end{cor}

The corollary is an immediate consequence of Theorem \ref{Thm:rep-ast-general} and the following lemma.

\begin{lem}
For any $k,r\in\mathbb{N}$, we have
\begin{align}
(r)^{\ast k}=&\sum\limits_{l=1}^k\sum\limits_{a_1+\cdots+a_l=k\atop a_i>0}\frac{k!}{a_1!\cdots a_l!}
(ra_1,ra_2,\ldots,ra_l),
\label{Eq:ast-r+}
\end{align}
and
\begin{align}
(\bar{r})^{\ast k}=&\sum\limits_{l=1}^k\sum\limits_{a_1+\cdots+a_l=k\atop a_i>0}\frac{k!}{a_1!\cdots a_l!}
(\widetilde{ra_1},\widetilde{ra_2},\ldots,\widetilde{ra_l}),
\label{Eq:ast-r-}
\end{align}
where $\widetilde{ra_i}$ is defined in \eqref{Eq:tilde}.
\end{lem}

\proof We give the proof of \eqref{Eq:ast-r-}. The proof of \eqref{Eq:ast-r+} is similar and easier. We use induction on $k$.
The initial value $k=1$ is trivial. We now assume that $k>0$. By induction assumption, we have
\begin{align*}
(\bar{r})^{\ast k}=(\bar{r})\ast \sum\limits_{l=1}^{k-1}\sum\limits_{a_1+\cdots+a_l=k-1\atop a_i>0}\frac{(k-1)!}{a_1!\cdots a_l!}
(\widetilde{ra_1},\widetilde{ra_2},\ldots,\widetilde{ra_l}).
\end{align*}
Now we have
\begin{align*}
(\bar{r})\ast (\widetilde{ra_1},\widetilde{ra_2},\ldots,\widetilde{ra_l})=&\sum\limits_{i=1}^{l+1} (\widetilde{ra_1},\ldots,\widetilde{ra_{i-1}},\bar{r},\widetilde{ra_{i+1}},\ldots,\widetilde{ra_l})\\
&+\sum\limits_{i=1}^l (\widetilde{ra_1},\ldots,\widetilde{ra_{i-1}},\widetilde{r(a_{i}+1)},\widetilde{ra_{i+1}},\ldots,\widetilde{ra_l}),
\end{align*}
which induces
\begin{align}
(\bar{r})^{\ast k}=&\sum\limits_{l=1}^{k}\sum\limits_{b_1+\cdots+b_l=k\atop \forall b_j>0, \exists b_i=1}\frac{(k-1)!}{b_1!\ldots b_l!}
(\widetilde{rb_1},\ldots,\widetilde{rb_l})
\label{Eq:justinLemma}\\
&+\sum\limits_{l=1}^k\sum\limits_{i=1}^l\sum\limits_{b_1+\cdots+b_l=k\atop
\forall b_j>0,b_i>1}\frac{(k-1)!b_i}{b_1!\cdots b_l!}(\widetilde{rb_1},\ldots,\widetilde{rb_l}).
\nonumber
\end{align}
The second term of the right-hand side of \eqref{Eq:justinLemma} equals
\begin{align*}
&\sum\limits_{l=1}^k\sum\limits_{i=1}^l\sum\limits_{b_1+\cdots+b_l=k\atop
\forall b_j>1}\frac{(k-1)!b_i}{b_1!\cdots b_l!}(\widetilde{rb_1},\ldots,\widetilde{rb_l})\\
&+\sum\limits_{l=1}^k\sum\limits_{a=1}^l\sum\limits_{i=1}^l\sum\limits_{b_1+\cdots+b_l=k\atop
\forall b_j>0,b_a=1,b_i>1}\frac{(k-1)!b_i}{b_1!\cdots b_l!}(\widetilde{rb_1},\ldots,\widetilde{rb_l}),
\end{align*}
which is
\begin{align*}
&\sum\limits_{l=1}^k\sum\limits_{b_1+\cdots+b_l=k\atop
\forall b_j>1}\frac{k!}{b_1!\cdots b_l!}(\widetilde{rb_1},\ldots,\widetilde{rb_l})\\
&+\sum\limits_{l=1}^k\sum\limits_{a=1}^l\sum\limits_{i=1\atop i\neq a}^l\sum\limits_{b_1+\cdots+b_l=k\atop
\forall b_j>0,b_a=1,b_i>1}\frac{(k-1)!b_i}{b_1!\cdots b_l!}(\widetilde{rb_1},\ldots,\widetilde{rb_l}).
\end{align*}
The equation \eqref{Eq:ast-r-} follows from the above formula and \eqref{Eq:justinLemma}.
\qed

We remark that one can also prove Corollary \ref{Cor:express-special} by using the multinomial theorem
$$(x_1+x_2+\cdots+x_n)^k=\sum\limits_{a_1+\cdots+a_n=k\atop a_i\geqslant 0}\frac{k!}{a_1!\cdots a_n!}x_1^{a_1}\cdots x_n^{a_n}.$$

We give two more examples. For $n\in\mathbb{Z}\backslash\{0\}$ with $n\geqslant 2$ if $n>0$, we have
\begin{align}
S(n;\{1\}^4)=&\zeta(n,4)+4\zeta(n,3,1)+4\zeta(n,1,3)+6\zeta(n,2,2)+12\zeta(n,2,1,1)\label{Eq:example-4}\\
&+12\zeta(n,1,2,1)+12\zeta(n,1,1,2)+24\zeta(n,1,1,1,1),\nonumber
\end{align}
and
\begin{align}
S(n;\{\bar{1}\}^4)=&\zeta(n,4)+4\zeta(n,\bar{3},\bar{1})+4\zeta(n,\bar{1},\bar{3})+6\zeta(n,2,2)+12\zeta(n,2,\bar{1},\bar{1})
\label{Eq:example-4-}\\
&+12\zeta(n,\bar{1},2,\bar{1})+12\zeta(n,\bar{1},\bar{1},2)+24\zeta(n,\bar{1},\bar{1},\bar{1},\bar{1}).\nonumber
\end{align}

In the end of this section, we state a result, which can be regarded as a generalization of \cite[Theorem 5.1 (i)]{Flajolet-Salvy}.

\begin{prop}\label{Prop:Combi-Zeta}
For any positive integers $n,m$ with $n\geqslant 2$, we have
\begin{align*}
\sum\limits_{l=1}^m\frac{(-1)^{m-l}}{l!}\sum\limits_{r_1+\cdots+r_l=m\atop r_i>0}\frac{1}{r_1\cdots r_l}
S(n;r_1,\ldots,r_l)=\zeta(n,\{1\}^m).
\end{align*}
\end{prop}

The above proposition is a direct consequence of Theorem \ref{Thm:rep-ast} and the following lemma.

\begin{lem}\label{Lem:ast}
For any $m\in\mathbb{N}$, we have
\begin{align}
\sum\limits_{l=1}^m\frac{(-1)^{m-l}}{l!}\sum\limits_{r_1+\cdots+r_l=m\atop r_i>0}\frac{1}{r_1\cdots r_l}
(r_1)\ast(r_2)\ast \cdots \ast(r_l)=(\{1\}^m),
\label{Eq:ast-sequences}
\end{align}
which holds in $\mathscr{S}^{+}\otimes_{\mathbb{Z}} \mathbb{Q}$, where $\mathscr{S}^{+}$ is the free
abelian group generated by all finite sequences of positive integers.
\end{lem}

\proof Let $\mathfrak{h}=\mathbb{Q}\langle x,y\rangle$ be the $\mathbb{Q}$-algebra of polynomials in two noncommutative variables $x$ and $y$. It has a subalgebra $\mathfrak{h}^1=\mathbb{Q}+\mathfrak{h}y$, which is a free algebra generated by $z_k=x^{k-1}y$ for $k=1,2,\ldots$. There is a commutative product $\ast$ in $\mathfrak{h}^1$, and we can identity $\mathscr{S}^{+}\otimes_{\mathbb{Z}}\mathbb{Q}$ with ($\mathfrak{h}^1,\ast$) by the identification $(k_1,\ldots,k_n)\mapsto z_{k_1}\cdots z_{k_n}$ (see \cite{Hoffman,IKZ}). Now the equivalent equation of \eqref{Eq:ast-sequences} in $\mathfrak{h}^1$ is the following
\begin{align}
\sum\limits_{l=1}^m\frac{(-1)^{m-l}}{l!}\sum\limits_{r_1+\cdots+r_l=m\atop r_i>0}\frac{1}{r_1\cdots r_l}
z_{r_1}\ast z_{r_2}\ast \cdots \ast z_{r_l}=y^m.
\label{Eq:ast-words}
\end{align}
We compute the generating function of the left-hand side of \eqref{Eq:ast-words} as the following
\begin{align*}
&1+\sum\limits_{m=1}^\infty\sum\limits_{l=1}^m\frac{(-1)^{m-l}}{l!}\sum\limits_{r_1+\cdots+r_l=m\atop r_i>0}\frac{1}{r_1\cdots r_l}
z_{r_1}\ast z_{r_2}\ast \cdots \ast z_{r_l}u^m\\
=&\sum\limits_{l=0}^\infty\frac{(-1)^l}{l!}\left(\sum\limits_{n=1}^\infty \frac{(-1)^{n}}{n}z_{n}u^{n}\right)^{\ast l}
=\exp_{\ast}\left(\sum\limits_{n=1}^\infty \frac{(-1)^{n-1}}{n}z_nu^n\right).
\end{align*}
Then equation \eqref{Eq:ast-words} is equivalent to
$$\exp_{\ast}\left(\sum\limits_{n=1}^\infty \frac{(-1)^{n-1}}{n}z_nu^n\right)=\frac{1}{1-yu},$$
which is a special case of \cite[Eq. 4.6]{IKZ}.
\qed

We list some examples. For $m=2$, we get
$$S(n;1,1)=\zeta(n,2)+2\zeta(n,1,1).$$
which is just the last equation in \cite{BBG}.
When $m=3$, we get
$$S(n;1,1,1)-3S(n;2,1)=-2\zeta(n,3)+6\zeta(n,1,1,1),$$
which can deduce the fact that if $n$ is even, then $S(n;1,1,1)-3S(n;2,1)$ is expressible in terms of zeta values
(\cite[Theorem 5.1(i)]{Flajolet-Salvy}).
When $m=4$, we get
$$S(n;\{1\}^4)-6S(n;2,1,1)+8S(n;3,1)+3S(n;2,2)=6\zeta(n,4)+24\zeta(n,\{1\}^4),$$
which implies that if $n$ is odd, then $S(n;\{1\}^4)-6S(n;2,1,1)+8S(n;3,1)+3S(n;2,2)$ is expressible in terms of zeta values.


\section{Some explicit evaluations}

We first give a typical example
\begin{align}
S(2;1,1)=\sum\limits_{n=1}^\infty\frac{\left(1+\frac{1}{2}+\cdots+\frac{1}{n}\right)^2}{(n+1)^2}=\frac{11}{4}\zeta(4),
\label{Eq:eva-example}
\end{align}
which is \cite[Eq.(2)]{BBG}. By \eqref{Eq:example-2}, we have
$$S(2;1,1)=\zeta(2,2)+2\zeta(2,1,1).$$
Hence \eqref{Eq:eva-example} follows from the facts $\zeta(2,2)=\frac{3}{4}\zeta(4)$ and $\zeta(2,1,1)=\zeta(4)$.
Through this example, we see that using Theorem \ref{Thm:rep-ast-general} to treat harmonic sums, one
should know the information of alternating Euler sums as much as possible.

\subsection{Relations and evaluations for multiple zeta values}

We first recall some relations among multiple zeta values, including duality formula, sum formula,
double shuffle relation and Aomoto-Drinfel'd-Zagier formula.

For a sequence
$\mathbb{\mathbf{k}}=(a_1+1,\{1\}^{b_1-1},a_2+1,\{1\}^{b_2-1},\ldots,a_{s}+1,\{1\}^{b_s-1})$
with $a_1,b_1,\ldots,a_s,b_s,s\in\mathbb{N}$, the dual sequence of $\mathbf{k}$ is defined as
$\mathbf{k}'=(b_s+1,\{1\}^{a_s-1},\ldots,b_2+1,\{1\}^{a_2-1},b_1+1,\{1\}^{a_1-1})$.
Then the duality formula (\cite{Zagier}) claims that
\begin{align}
\zeta(\mathbf{k})=\zeta(\mathbf{k}').
\label{Eq:duality}
\end{align}

Sum formula was conjectured in \cite{Zagier} and first proved in \cite{Granville},
which says that the sum of multiple zeta values of fixed weight and depth is
independent of depth. More explicitly, for any $k,n\in\mathbb{N}$ with $k>n$, we have
\begin{align}
\sum\limits_{k_1+\cdots+k_n=k\atop k_1\geqslant 2,k_2,\ldots,k_n\geqslant 1}\zeta(k_1,\ldots,k_n)=\zeta(k).
\label{Eq:sum}
\end{align}

Double shuffle relation is deduced from the fact that we have two ways to write a product of two multiple zeta values
as a $\mathbb{Z}$-linear combination of multiple zeta values. We may state this relation as
\begin{align}
\zeta(\mathbf{k})\zeta(\mathbf{l})=\zeta(\mathbf{k}\ast\mathbf{l})=\zeta(\mathbf{k}\shuffle \mathbf{l}),
\label{Eq:doubleshuffle}
\end{align}
where $\mathbf{k}=(k_1,\ldots,k_n),\mathbf{l}=(l_1,\ldots,l_m)$ are sequences of positive integers with $k_1,l_1\geqslant 2$,
and $\shuffle$ is the shuffle product. For details, one can refer to \cite{IKZ}.

The Aomoto-Drinfel'd-Zagier formula reads
\begin{align}
\sum\limits_{m,n=1}^\infty \zeta(m+1,\{1\}^{n-1})x^my^n=1-\exp\left(\sum\limits_{n=2}^\infty\zeta(n)\frac{x^n+y^n-(x+y)^n}{n}\right),
\label{Eq:ADZ-formula}
\end{align}
which implies that for any $m,n\in\mathbb{N}$, the multiple zeta value $\zeta(m+1,\{1\}^{n-1})$ can be represented as a polynomial of zeta values with rational coefficients. In particular, one can find explicit formulas for small weights.

Now we list some evaluations of multiple zeta values via zeta values, which will be used in the next subsection.
For weight $5$, we have
\begin{align}
&\zeta(4,1)=-\zeta(2)\zeta(3)+2\zeta(5),
\label{Eq:formula-z41}\\
&\zeta(3,2)=3\zeta(2)\zeta(3)-\frac{11}{2}\zeta(5),
\label{Eq:formula-z32}\\
&\zeta(2,3)=-2\zeta(2)\zeta(3)+\frac{9}{2}\zeta(5).
\label{Eq:formula-z23}
\end{align}

We give the proof of the above evaluations. The evaluation \eqref{Eq:formula-z41} follows from Aomoto-Drinfel'd-Zagier formula. Double shuffle relation
$\zeta(2)\zeta(3)=\zeta(2,3)+\zeta(3,2)+\zeta(5)=\zeta(2,3)+3\zeta(3,2)+6\zeta(4,1)$
and \eqref{Eq:formula-z41} deduce \eqref{Eq:formula-z32} and \eqref{Eq:formula-z23}.

For weight $6$, we have
\begin{align}
&\zeta(5,1)=-\frac{1}{2}\zeta(3)^2+\frac{3}{4}\zeta(6),
\label{Eq:formula-z51}\\
&\zeta(4,2)=\zeta(3)^2-\frac{4}{3}\zeta(6),
\label{Eq:formula-z42}\\
&\zeta(3,3)=\frac{1}{2}\zeta(3)^2-\frac{1}{2}\zeta(6),
\label{Eq:formula-z33}\\
&\zeta(2,4)=-\zeta(3)^2+\frac{25}{12}\zeta(6),
\label{Eq:formula-z24}\\
&\zeta(4,1,1)=-\zeta(3)^2+\frac{23}{16}\zeta(6),
\label{Eq:formula-z411}\\
&\zeta(3,2,1)=3\zeta(3)^2-\frac{203}{48}\zeta(6),
\label{Eq:formula-z321}\\
&\zeta(3,1,2)=-\frac{3}{2}\zeta(3)^2+\frac{53}{24}\zeta(6),
\label{Eq:formula-z312}\\
&\zeta(2,2,2)=\frac{3}{16}\zeta(6).
\label{Eq:formula-z222}
\end{align}

The evaluations  \eqref{Eq:formula-z51}
and \eqref{Eq:formula-z411} follow from Aomoto-Drinfel'd-Zagier formula under the helps of the facts
$\zeta(2)=\pi^2/6$, $\zeta(4)=\pi^4/90$ and $\zeta(6)=\pi^6/945$. The harmonic shuffle product
$\zeta(3)\zeta(3)=2\zeta(3,3)+\zeta(6)$
gives \eqref{Eq:formula-z33}. The double shuffle relation
$\zeta(2)\zeta(4)=\zeta(2,4)+\zeta(4,2)+\zeta(6)=\zeta(2,4)+4\zeta(4,2)+2\zeta(3,3)+8\zeta(5,1),$
and evaluations for $\zeta(3,3),\zeta(5,1),\zeta(2),\zeta(4),\zeta(6)$ deduce \eqref{Eq:formula-z42} and
\eqref{Eq:formula-z24}. Using the duality $\zeta(2,3,1)=\zeta(3,1,2)$, we have the double shuffle relation
\begin{align*}
\zeta(2)\zeta(3,1)=&\zeta(3,2,1)+2\zeta(3,1,2)+\zeta(5,1)+\zeta(3,3)\\
=&4\zeta(3,2,1)+2\zeta(3,1,2)+9\zeta(4,1,1),
\end{align*}
which together with evaluations for $\zeta(5,1),\zeta(4,1,1),\zeta(2),\zeta(3,1)$ imply \eqref{Eq:formula-z321} and \eqref{Eq:formula-z312}.
Finally, the evaluation \eqref{Eq:formula-z222} is a special case of the following more general formula
$$\zeta(\{2\}^n)=\frac{\pi^{2n}}{(2n+1)!},\;(n\in\mathbb{N}).$$

The above evaluations should have appeared in the existing literatures.
We list and prove them here is to emphasize that we can obtain them formally.

\subsection{Evaluations for $S(n;\{1\}^k)$}

We give the formal proofs of some evaluations of $S(n;\{1\}^k)$ listed in \cite[Table 3]{BaileyBG} of weights up to $6$.

\begin{prop}
The following evaluations hold
\begin{align}
S(2;\{1\}^3)=&\zeta(2)\zeta(3)+\frac{15}{2}\zeta(5),\label{Eq:formula-s23}\\
S(3;\{1\}^3)=&-\frac{33}{16}\zeta(6)+2\zeta(3)^2,\label{Eq:formula-s33}\\
S(2;\{1\}^4)=&3\zeta(3)^2+\frac{859}{24}\zeta(6).\label{Eq:formula-s24}
\end{align}
\end{prop}
Note that \eqref{Eq:formula-s24} is an experimental detected formula  in \cite{BaileyBG}.

\proof
We have
$$S(2;\{1\}^3)=\zeta(2,3)+3\zeta(2,2,1)+3\zeta(2,1,2)+6\zeta(2,1,1,1).$$
Duality formula gives that
$\zeta(2,2,1)=\zeta(3,2)$, $\zeta(2,1,2)=\zeta(2,3)$ and $\zeta(2,1,1,1)=\zeta(5)$.
Hence we have
$$S(2;\{1\}^3)=4\zeta(2,3)+3\zeta(3,2)+6\zeta(5).$$
Now we get \eqref{Eq:formula-s23} by \eqref{Eq:formula-z32} and \eqref{Eq:formula-z23}.

We have
$$S(3;\{1\}^3)=\zeta(3,3)+3\zeta(3,2,1)+3\zeta(3,1,2)+6\zeta(3,1,1,1).$$
Then the duality $\zeta(3,1,1,1)=\zeta(5,1)$ and \eqref{Eq:formula-z51},
\eqref{Eq:formula-z321}, \eqref{Eq:formula-z312} imply \eqref{Eq:formula-s33}.

Finally, we have
\begin{align*}
S(2;\{1\}^4)=&\zeta(2,4)+4\zeta(2,3,1)+4\zeta(2,1,3)+6\zeta(2,2,2)+12\zeta(2,2,1,1)\\
&+12\zeta(2,1,2,1)+12\zeta(2,1,1,2)+24\zeta(2,1,1,1,1).
\end{align*}
Sum formula gives
$$\zeta(2,2,1,1)+\zeta(2,1,2,1)+\zeta(2,1,1,2)=\zeta(6)-\zeta(3,1,1,1),$$
and
$$\zeta(2,3,1)+\zeta(2,1,3)+\zeta(2,2,2)=\zeta(6)-\zeta(3,2,1)-\zeta(3,1,2)-\zeta(4,1,1).$$
Hence we have
\begin{align*}
S(2;\{1\}^4)=&40\zeta(6)-12\zeta(5,1)+\zeta(2,4)-4\zeta(4,1,1)-4\zeta(3,2,1)\\
&-4\zeta(3,1,2)+2\zeta(2,2,2).
\end{align*}
Here we have used the dualities $\zeta(3,1,1,1)=\zeta(5,1)$ and $\zeta(2,1,1,1,1)=\zeta(6)$.
Now the evaluation \eqref{Eq:formula-s24} follows from \eqref{Eq:formula-z51}, \eqref{Eq:formula-z24}, \eqref{Eq:formula-z411},
\eqref{Eq:formula-z321}, \eqref{Eq:formula-z312} and \eqref{Eq:formula-z222}.
\qed

Of course, we think one can continue to give formal proofs of the other experimentally detected formulas for $S(n;\{1\}^k)$ in
\cite[Table 3]{BaileyBG}. In fact, we have checked that the formulas for $S(4;\{1\}^3)$, $S(3,\{1\}^4)$ and $S(2;\{1\}^5)$ in \cite[Table 3]{BaileyBG} should be proved if we adopt the evaluation formulas for triple zeta values of weight $7$ given in \cite{BG}.

\subsection{Alternating harmonic sums}

We present here, as an example, the evaluations of alternating harmonic sums of weight $3$.

As done in above, we have to obtain the evaluations of alternating Euler sums of weight $3$.
First, note that $\zeta(\bar{1})=-\log 2$ and $\zeta(\bar{n})=(2^{1-n}-1)\zeta(n)$ for $n\in\mathbb{N}$ with $n\geqslant 2$.
Then using (regularized) double shuffle relation(\cite{BJOP,Zhao}), we obtain the following formulas:
\begin{align*}
\zeta(\bar{1},1)=&\frac{1}{2}\log^2 2,&
\zeta(\bar{1},\bar{1})=&\frac{1}{2}\log^2 2-\frac{1}{2}\zeta(2),
\end{align*}
and
\begin{align*}
&\zeta(2,\bar{1})=-\frac{3}{2}\zeta(2)\log 2+\zeta(3),&
&\zeta(\bar{2},1)=\frac{1}{8}\zeta(3),\\
&\zeta(\bar{2},\bar{1})=\frac{3}{2}\zeta(2)\log 2-\frac{13}{8}\zeta(3),&
&\zeta(\bar{1},2)=\frac{1}{2}\zeta(2)\log 2-\frac{1}{4}\zeta(3),\\
&\zeta(\bar{1},\bar{2})=-\zeta(2)\log 2+\frac{5}{8}\zeta(3),&
&\zeta(\bar{1},1,1)=-\frac{1}{6}\log^3 2,\\
&\zeta(\bar{1},1,\bar{1})=-\frac{1}{6}\log^3 2+\frac{1}{8}\zeta(3),\\
&\zeta(\bar{1},\bar{1},1)=-\frac{1}{6}\log^3 2+\frac{1}{2}\zeta(2)\log 2-\frac{7}{8}\zeta(3),\\
&\zeta(\bar{1},\bar{1},\bar{1})=-\frac{1}{6}\log^3 2+\frac{1}{2}\zeta(2)\log 2-\frac{1}{4}\zeta(3).
\end{align*}
Note that they are the same as that given in \cite{BJOP}.

The two formulas of weight $2$ follow from the double shuffle relation
$\zeta(\bar{1})\zeta(\bar{1})=2\zeta(\bar{1},\bar{1})+\zeta(2)=2\zeta(\bar{1},1)$. For the weight $3$ cases, we have
$4$ pairs double shuffle relations: $\zeta(2)\zeta(\bar{1})=\zeta(2,\bar{1})+\zeta(\bar{1},2)+\zeta(\bar{3})= \zeta(2,\bar{1})+\zeta(\bar{2},\bar{1})+\zeta(\bar{1},\bar{2})$, $\zeta(\bar{2})\zeta(\bar{1})
=\zeta(\bar{2},\bar{1})+\zeta(\bar{1},\bar{2})+\zeta(3)=\zeta(\bar{1},2)+2\zeta(\bar{2},1)$,
$\zeta(\bar{1})\zeta(\bar{1},1)=2\zeta(\bar{1},\bar{1},1)+\zeta(\bar{1},1,\bar{1})+\zeta(2,1)+\zeta(\bar{1},\bar{2})
=3\zeta(\bar{1},1,1)$, $\zeta(\bar{1})\zeta(\bar{1},\bar{1})=3\zeta(\bar{1},\bar{1},\bar{1})+\zeta(2,\bar{1})+\zeta(\bar{1},2)
=2\zeta(\bar{1},1,\bar{1})+\zeta(\bar{1},\bar{1},\bar{1})$. By computing the shuffle and harmonic shuffle of $(1)$ and $(\bar{2})$,
we get a regularized double shuffle relation $\zeta(\bar{2},1)+\zeta(\bar{3})=\zeta(2,\bar{1})+\zeta(\bar{2},\bar{1})$. Then using
these $9$ equations and the evaluations for weight $2$, we obtain the evaluations for weight $3$.

Now we have the evaluations for alternating harmonic sums of weight $3$.

\begin{prop}
We have the evaluations
\begin{align*}
&S(2;\bar{1})=\frac{3}{2}\zeta(2)\log 2-\zeta(3),\\
&S(\bar{2};1)=\frac{1}{8}\zeta(3),\\
&S(\bar{2};\bar{1})=-\frac{3}{2}\zeta(2)\log 2+\frac{13}{8}\zeta(3),\\
&S(\bar{1};2)=\frac{1}{2}\zeta(2)\log 2-\frac{1}{4}\zeta(3),\\
&S(\bar{1};1,1)=-\frac{1}{3}\log^3 2+\frac{1}{2}\zeta(2)\log 2-\frac{1}{4}\zeta(3),\\
&S(\bar{1};1,\bar{1})=\frac{1}{3}\log^3 2+\frac{1}{2}\zeta(2)\log 2+\frac{1}{8}\zeta(3),\\
&S(\bar{1};\bar{1},\bar{1})=-\frac{1}{3}\log^3 2+\frac{3}{2}\zeta(2)\log 2-\frac{3}{4}\zeta(3).
\end{align*}
\end{prop}

\proof The evaluations follow from Theorem \ref{Thm:rep-ast-general} and evaluation formulas for alternating Euler sums of weight $3$ given above.
\qed

We finally remark that using the table in \cite{BJOP} for alternating Euler sums of weight $4$, one can obtain evaluation formulas
for harmonic sums of weight $4$ in terms of $\log 2,\zeta(2),\zeta(3),\zeta(4)$ and an undetermined sum $\zeta(\bar{3},1)$,
which can be represented via $\Li_{4}(1/2)$ by using duality formula for $\zeta(\bar{1},\bar{1},1,1)$. 

$ $

\vskip5pt





\begin{thebibliography}{99}

\bibitem{BaileyBG} D. H. Bailey, J. M. Borwein and R. Girgensohn, Experimental evaluation of Euler sums, Experimental Math. 3(1)(1994), 17-30.

\bibitem{BJOP} M. Bigotte, G. Jacob, N. E. Oussous and M. Petitot, Lyndon words and shuffle algebras for generating the coloured multiple zeta
values relations tables, Theoretical Computer Sciences 273(2002), 271-282.

\bibitem{BB} D. Borwein and J. M. Borwein, On an intriguing integral and some series related to $\zeta(4)$, Proc. Amer. Math. Soc. 123(4)(1995), 1191-1198.

\bibitem{BBG} D. Borwein, J. M. Borwein and R. Girgensohn, Explicit evaluation of Euler sums, Proc. Edinburgh Math. Soc. 38(1995), 277-294.

\bibitem{BG} J. M. Borwein and R. Girgensohn, Evaluation of triple Euler sums, Electron. J. Combin. 3(1996), \#R23, 1-27.

\bibitem{BorweinBradleyB} J. M. Borwein, D. M. Bradley and D. J. Broadhurst, Evaluations of $k$-fold Euler/Zagier sums: a compendium of results
for arbitrary $k$, Electron. J. Combin. 4(2)(1997), \#R5, 1-21.

\bibitem{Choi-Srivastava} J. Choi and H. M. Srivastava, Explicit evaluation of Euler and related sums, Ramanujan J. 10(2005), 51-70.

\bibitem{Chu} W. Chu, Hypergeometric series and the Riemann zeta function, Acta Arith. 82(1997), 103-118.

\bibitem{Doelder} P. J. de Doelder, On some series containing $\psi(x)-\psi(y)$ and $(\psi(x)-\psi(y))^2$ for certain values of $x$ and $y$,
J. Comput. Appl. Math. 37(1991), 125-141.

\bibitem{Flajolet-Salvy} P. Flajolet and B. Salvy, Euler sums and contour integral representations, Experimental Math. 7(1)(1998), 15-35.

\bibitem{Granville} A. Granville, A decomposition of Riemann zeta-function, in Analytic Number Theory,
London Math. Soc. Lecture Notes Ser., Vol. 247, Cambridge University Press,
Cambridge, 1997, 95-101.

\bibitem{Hoffman} M. E. Hoffman, The algebra of multiple harmonic series, Journal of Algebra,
194(1997), 477-495.

\bibitem{IKZ} K. Ihara, M. Kaneko and D. Zagier, Derivation and double shuffle relations for multiple zeta values, Compositio Math. 142(2006),
307-338.

\bibitem{Shen} L. -C. Shen, Remarks on some integrals and series involving the Stirling numbers and $\zeta(n)$, Trans. Amer. Math. Soc. 347(4)(1995), 1391-1399.

\bibitem{Zagier} D. Zagier, Values of zeta functions and their applications,
in First European Congress of Mathematics, Vol. II (Paris, 1992), volume 120
of Progr. Math., pages 497-512. Birkh\"{a}user, Basel, 1994.

\bibitem{Zhao} J. Zhao, On a conjecture of Borwein, Bradley and Broadhurst, J. reine angew. Math. 639(2010), 223-233.

\end{thebibliography}
\end{document}